\documentclass[11pt]{article}

\usepackage[a4paper,margin=1in]{geometry}
\usepackage[T1]{fontenc}
\usepackage{lmodern}
\usepackage{microtype}
\usepackage{amsmath,amssymb,amsthm,mathtools}
\usepackage{mathrsfs}
\usepackage{enumitem}
\usepackage{xcolor}
\usepackage{hyperref}
\usepackage[nameinlink,noabbrev]{cleveref}

\hypersetup{
  colorlinks=true,
  linkcolor=blue!60!black,
  citecolor=blue!60!black,
  urlcolor=blue!60!black,
  pdftitle={Degree of Maps and Stein Degree},
  pdfauthor={},
}

\theoremstyle{definition}
\newtheorem{definition}{Definition}[section]
\newtheorem{example}[definition]{Example}
\newtheorem{remark}[definition]{Remark}
\newtheorem{conjecture}[definition]{Conjecture}

\theoremstyle{plain}
\newtheorem{theorem}[definition]{Theorem}
\newtheorem{fact}[definition]{Fact}

\newcommand{\A}{\mathbb{A}}
\newcommand{\Pp}{\mathbb{P}}
\newcommand{\Q}{\mathbb{Q}}
\newcommand{\C}{\mathbb{C}}
\newcommand{\F}{\mathbb{F}}
\newcommand{\Spec}{\operatorname{Spec}}
\newcommand{\Nklt}{\operatorname{Nklt}}
\newcommand{\sdeg}{\operatorname{sdeg}}

\newcommand{\id}{\mathrm{id}}
\newcommand{\OO}{\mathcal{O}}

\title{\textbf{Stein degree of proper morphisms}}
\author{Caucher Birkar}
\date{}

\begin{document}

\maketitle

\begin{abstract}
The notion of degree begins in field theory as the dimension of a field extension. In algebraic geometry, this idea reappears as the degree of a finite morphism, defined using the induced extension of function fields. For proper morphisms that are not necessarily finite, Stein factorization isolates the finite part of the map and leads to the notion of \emph{Stein degree}. This invariant is especially useful in birational geometry, where it interacts naturally with singularities of pairs and the study of log Calabi--Yau fibrations. In this article we give an expository introduction to these ideas, discuss motivating examples, and explain a boundedness problem for Stein degree arising in recent work of the author and collaborators.
\end{abstract}

\tableofcontents

\section{Introduction}

This article is an expanded expository version of a general seminar on degree of maps and Stein degree presented at the International Congress of Chinese Mathematicians in January 2026 in Shanghai, China.

One of the most basic numerical invariants in algebra is the degree of a field extension. If a field \(L\) contains a field \(K\), then \(L\) is naturally a vector space over \(K\), and its dimension measures how much larger \(L\) is than \(K\). In algebraic geometry the same idea appears in a geometric form: if \(X\to Y\) is a finite morphism of algebraic varieties, then the induced extension of function fields \(K(Y)\subseteq K(X)\) is finite, and its degree gives the degree of the morphism.

Many morphisms of interest, however, are proper rather than finite. In this setting the correct replacement for ordinary degree is obtained by using Stein factorization, which decomposes a proper morphism into a map with connected fibers followed by a finite map. The degree of that finite part is called the \emph{Stein degree}, introduced in \cite{B-moduli}. This is a basic invariant and it will not be suprising if it appears in many areas of algebraic geometry.

The goal of this article is to explain the passage
\[
\text{field degree} \longrightarrow \text{degree of finite morphisms} \longrightarrow \text{Stein degree of proper morphisms},
\]
and then to describe how Stein degree appears in birational geometry, especially in the study of singularities of pairs and log Calabi--Yau fibrations.

Stein degree is particularly useful in birational geometry, where it interacts with singularities of pairs and boundedness problems. The resulting questions are interesting not only geometrically but also arithmetically, especially when the relevant loci are defined over non-algebraically closed fields. 

The boundedness of Stein degree is therefore a natural and subtle problem. It lies at the intersection of several major themes of modern algebraic geometry: the minimal model program, Fano and Calabi-Yau varieties and complements, toroidal methods, moduli theory, arithmetic geometry, and motivic invariants. It is rapidly developing into a fruitful research direction \cite{B-moduli, BQ-sdeg, BQ-sdeg-non-Fano, LS-motivic}.

\section{Degree of field extensions}

We begin with the algebraic notion from which everything else is derived.

\begin{definition}
Let \(K\subseteq L\) be a field extension. The \emph{degree} of the extension is
\[
[L:K] := \dim_K L.
\]
If this dimension is finite, then the extension is called \emph{finite}.
\end{definition}

\begin{example}
Consider the extension
\[
\Q \subseteq \Q(\sqrt{2}).
\]
A basis of \(\Q(\sqrt{2})\) over \(\Q\) is given by \(\{1,\sqrt{2}\}\). Hence
\[
[\Q(\sqrt{2}):\Q]=2.
\]
\end{example}

\begin{example}
Similarly,
\[
\Q \subseteq \Q(\sqrt[3]{2})
\]
has basis \(\{1,\sqrt[3]{2},\sqrt[3]{4}\}\), so
\[
[\Q(\sqrt[3]{2}):\Q]=3.
\]
\end{example}

\begin{example}
For a prime \(p\), the finite field extension
\[
\F_p \subseteq \F_{p^n}
\]
has degree \(n\).
\end{example}

\begin{example}
Not every field extension is finite. For instance,
\[
[\C:\Q]=\infty.
\]
\end{example}

These examples show that degree measures algebraic complexity. In geometry, the analogous quantity is extracted from function fields.

\section{Finite morphisms of varieties}

Let \(f\colon X\to Y\) be a morphism of varieties. We recall the basic definition.

\begin{definition}
A morphism \(f\colon X\to Y\) is \emph{finite} if for every affine open subset \(V\subseteq Y\), the inverse image \(f^{-1}(V)\) is affine and
\[
\OO_X(f^{-1}(V))
\]
is a finite \(\OO_Y(V)\)-module.
\end{definition}

This is the geometric analogue of a finite ring extension. Finite morphisms are proper and have finite fibers, so they are natural geometric generalizations of finite field extensions.

A finite morphism \(f\colon X\to Y\) of vareities of the same dimension induces an inclusion of function fields
\[
K(Y)\subseteq K(X).
\]
The extension is finite, and its degree is the right notion of degree for the map.

\begin{definition}
Let \(f\colon X\to Y\) be a finite morphism of varieties with \(\dim X=\dim Y\). The \emph{degree} of \(f\) is
\[
\deg(f):=[K(X):K(Y)].
\]
\end{definition}

For a general point of \(Y\), this degree can be interpreted as the number of points in the fiber, counted appropriately.

\begin{example}
Consider
\[
f\colon \A^1\to \A^1,\qquad t\mapsto t^2.
\]
On function fields we have the inclusion
\[
k(t^2)\subseteq k(t).
\]
Since \(t\) satisfies the polynomial \(x^2-t^2=0\) over \(k(t^2)\), the extension has degree \(2\), and therefore
\[
\deg(f)=[k(t):k(t^2)]=2.
\]
Geometrically, a general point has two preimages.
\end{example}

\begin{example}
Now consider
\[
f\colon \A^1\to \A^1,\qquad t\mapsto t^3+t.
\]
Then
\[
k(t^3+t)\subseteq k(t),
\]
and \(t\) is algebraic of degree \(3\) over \(k(t^3+t)\). Hence
\[
\deg(f)=[k(t):k(t^3+t)]=3.
\]
\end{example}

\section{Further examples of finite morphisms}

The definition becomes more geometric when applied to singular varieties and projections.

\begin{example}[Normalization of a cusp]
Consider the cusp
\[
V(y^2-x^3)\subset \A^2
\]
and the map
\[
f\colon \A^1\to V(y^2-x^3),\qquad t\mapsto (t^2,t^3).
\]
This is the normalization map.

On function fields,
\[
k(t^2,t^3)=k(t),
\]
because \(t=t^3/t^2\) in the function field. Thus
\[
[k(t):k(t^2,t^3)]=1,
\]
so
\[
\deg(f)=1.
\]
This shows that a finite map of degree \(1\) need not be an isomorphism; rather, it is birational.
\end{example}

\begin{example}[Projection from a plane curve]
Let
\[
X=V(y^2-x^3+x)\subset \A^2,
\]
and consider the projection
\[
f\colon X\to \A^1,\qquad (x,y)\mapsto x.
\]
The defining equation gives
\[
y^2=x^3-x,
\]
so over the rational function field \(k(x)\), the coordinate \(y\) satisfies a quadratic equation. Therefore
\[
[K(X):k(x)]=2,
\]
and hence
\[
\deg(f)=2.
\]
\end{example}

\begin{example}[Power map on affine space]
Consider
\[
f\colon \A^n\to \A^n,\qquad (x_1,\dots,x_n)\mapsto (x_1^5,\dots,x_n^5).
\]
If the target coordinates are \(y_i=x_i^5\), then
\[
k(y_1,\dots,y_n)\subseteq k(x_1,\dots,x_n),
\]
and each variable contributes degree \(5\). Hence
\[
[k(x_1,\dots,x_n):k(y_1,\dots,y_n)]=5^n,
\]
so
\[
\deg(f)=5^n.
\]
\end{example}

\section{Stein factorization}

The degree of a finite morphism is very natural, but proper morphisms need not be finite. Stein factorization extracts the finite part of a proper map.

\begin{fact}[Stein factorization]
Let \(f\colon X\to Y\) be a proper morphism of Noetherian schemes. Then there exists a unique factorization
\[
X \xrightarrow{g} Z \xrightarrow{h} Y
\]
such that:
\begin{enumerate}[label=\textup{(\roman*)}]
    \item \(g\) is proper with connected fibers,
    \item \(h\) is finite,
    \item \(\OO_Z \cong g_*\OO_X\).
\end{enumerate}
\end{fact}

The idea is that the map \(g\) contracts each connected component of a fiber of \(f\) to a single point, while \(h\) records the remaining finite behavior over the base.

\begin{example}
Consider
\[
X=V(x^2+y^2+tz^2)\subset \Pp^2\times \A^1
\]
and the morphism
\[
f\colon X\to \A^1,\qquad ([x:y:z],t)\mapsto t^2.
\]
The projection \(X\to \A^1_t\) has connected fibers, while the map \(t\mapsto t^2\) is finite of degree \(2\). Thus the Stein factorization of \(f\) is
\[
X\xrightarrow{g}\A^1\xrightarrow{h}\A^1,
\]
where \(g\) is projection to the \(t\)-coordinate and
\[
h(t)=t^2.
\]
\end{example}

This example is simple but very instructive: the proper morphism \(f\) is not finite, yet its ``finite content'' is exactly the degree-two map on the base.

\section{Stein degree}

Motivated by Stein factorization, one defines a degree for proper morphisms.

\begin{definition}
Let \(f\colon X\to Y\) be a proper morphism of varieties over a field \(k\), and let
\[
X\xrightarrow{g}Z\xrightarrow{h}Y
\]
be its Stein factorization. The \emph{Stein degree} of \(X\) over \(Y\) is
\[
\sdeg(X/Y)=
\begin{cases}
\deg(h) & \text{if } f \text{ is surjective},\\[4pt]
0 & \text{otherwise}.
\end{cases}
\]
\end{definition}

Thus Stein degree reduces the study of a proper map to the ordinary degree of the finite morphism appearing in its Stein factorization.

\begin{remark}
This invariant was introduced in \cite{B-moduli} in the construction of moduli of stable minimal models. It is useful precisely because many naturally occurring morphisms in birational geometry are proper but not finite.
\end{remark}

\begin{example}[Closed points]
Let \(x\in X\) be a closed point and let \(y=f(x)\in Y\). Then
\[
\sdeg(x/y)=[k(x):k(y)].
\]
So Stein degree recovers the arithmetic degree of the residue field extension.
\end{example}

\begin{example}[Birational proper morphisms]
If \(f\colon X\to Y\) is a proper birational morphism, then the finite map in the Stein factorization has degree \(1\). Hence
\[
\sdeg(X/Y)=1.
\]
\end{example}

\section{Examples of Stein degree}

\begin{example}
Consider
\[
f\colon \Pp^1\to \Pp^1,\qquad [x:y]\mapsto [x^2:y^2].
\]
This map is finite and surjective of degree \(2\). Its Stein factorization is simply
\[
\Pp^1 \xrightarrow{\id} \Pp^1 \xrightarrow{h} \Pp^1,
\]
where \(h\) is the same degree-two map. Therefore
\[
\sdeg(\Pp^1/\Pp^1)=2.
\]
\end{example}

\begin{example}
Let
\[
X=V(x^2+y^2+tz^2)\subset \Pp^2\times \A^1,
\]
and define
\[
f([x:y:z],t)=t^n.
\]
Then the Stein factorization has finite part
\[
h\colon \A^1\to \A^1,\qquad t\mapsto t^n.
\]
Thus
\[
\sdeg(X/\A^1)=n.
\]
\end{example}

These examples show that Stein degree measures the residual finite covering behavior of a proper morphism after connected fibers have been collapsed.

\section{Singularities of pairs}

We now move to the birational-geometric setting in which Stein degree becomes especially interesting.

\begin{definition}
A \emph{pair} \((X,B)\) consists of:
\begin{itemize}
    \item a normal variety \(X\), and
    \item a \(\Q\)-divisor
    \[
    B=\sum a_i B_i,\qquad a_i\in \Q,
    \]
    such that \(K_X+B\) is \(\Q\)-Cartier.
\end{itemize}
\end{definition}

To measure the singularities of a pair, one pulls back to a log resolution.

\begin{definition}
Let \(f\colon W\to X\) be a log resolution of a pair \((X,B)\), and write
\[
K_W+B_W=f^*(K_X+B).
\]
Then:
\begin{itemize}
    \item \((X,B)\) is \emph{klt} if all coefficients of \(B_W\) are \(<1\),
    \item \((X,B)\) is \emph{lc} if all coefficients of \(B_W\) are \(\le 1\).
\end{itemize}
\end{definition}

The non-klt locus records where the singularities fail to be klt.

\begin{definition}
The \emph{non-klt locus} of a pair \((X,B)\) is
\[
\Nklt(X,B)=\{x\in X : (X,B)\text{ is not klt at }x\}.
\]
Equivalently,
\[
\Nklt(X,B)=f(\lfloor B_W\rfloor).
\]
\end{definition}

\subsection*{Basic examples}

\begin{example}
Let \(X\) be a smooth curve and let
\[
B=\frac12 p_1+p_2,
\]
where \(p_1,p_2\) are distinct points. Then \((X,B)\) is lc, but not klt at \(p_2\). Therefore
\[
\Nklt(X,B)=\{p_2\}.
\]
\end{example}

\begin{example}
Let \(X=\Pp^2\), and let
\[
B=\frac23L_1+\frac23L_2+\frac23L_3,
\]
where the lines \(L_i\) are distinct and all pass through a point \(p\). If one blows up \(p\), the exceptional divisor \(E\) appears with coefficient
\[
\mu_E B_W=1.
\]
Hence \((X,B)\) is lc but not klt at \(p\), and
\[
\Nklt(X,B)=\{p\}.
\]
\end{example}

\begin{example}
Let \(X=\Pp^2\) and let \(B=2L\), where \(L\) is a line. Then the pair is not klt along the whole line, so
\[
\Nklt(X,B)=L.
\]
\end{example}

\begin{example}
Some standard examples from the literature are:
\begin{itemize}
    \item the cone over a rational curve is klt,
    \item the cone over an elliptic curve is lc but not klt,
    \item quotient singularities are klt,
    \item toric pairs are lc,
    \item Du Val surface singularities are klt.
\end{itemize}
\end{example}

\section{The connectedness principle}

A fundamental theorem in the subject is the Shokurov--Koll\'ar connectedness theorem.

\begin{theorem}[Connectedness principle]
Let \((X,B)\) be a pair and let \(f\colon X\to Y\) be a proper morphism with connected fibers. If
\[
-(K_X+B)
\]
is \(f\)-ample, then \(\Nklt(X,B)\) is connected in a neighborhood of each fiber of \(f\).
\end{theorem}

This theorem is one of the main reasons the non-klt locus is such a useful object. It says that under positivity assumptions, the bad singularities cannot split into several disconnected pieces over a fiber.

From the perspective of Stein degree, one may think of this as saying that the relevant non-klt locus has Stein degree \(1\) over its image (but one needs to be careful as the non-klt locus may not be irreducible).

\begin{example}
Let
\[
X=\Pp^2,\qquad B=\frac23L_1+\frac23L_2+\frac23L_3,
\]
where the three lines pass through the same point \(p\). Since \(K_X=-3H\), we have
\[
-(K_X+B)=3H-\frac23(L_1+L_2+L_3)\equiv H,
\]
which is ample. As above,
\[
\Nklt(X,B)=L_1\cap L_2\cap L_3=\{p\},
\]
and this set is connected.
\end{example}

\begin{example}
Let
\[
X=\Pp^2,\qquad B=\frac23(C_1+C_2+C_3)
\]
where \(C_i\) are general curves of degree $\ge 2$ passing through \(n\ge 2\) fixed points. Then
\[
-(K_X+B)
\]
is not ample, and
\[
\Nklt(X,B)=\bigcap C_i
\]
consists of \(n\) points. In particular, the non-klt locus is not connected. This shows that the ampleness assumption in the theorem is essential.
\end{example}

\section{Log Calabi--Yau fibrations and Stein degree}

We now come to the boundedness questions that motivate this discussion.

\begin{definition}
A \emph{log Calabi--Yau fibration} (abbreviated lcyf) is a projective morphism
\[
(X,B)\to Y
\]
such that:
\begin{itemize}
    \item \((X,B)\) is lc,
    \item \(X\to Y\) is a contraction, i.e. $f_*\mathcal{O}_X=Y$ so it has connected fibers,
    \item \(K_X+B\sim_\Q 0/Y\).
\end{itemize}
\end{definition}

Suppose \(S\subset B\) is a component with coefficient \(\mu_S B=1\), and assume that \(S\) is \emph{horizontal} over \(Y\), meaning that it dominates \(Y\). One can then ask whether the Stein degree
\[
\sdeg(S/Y)
\]
is bounded in terms of the dimension.

\begin{theorem}[Birkar, \cite{B-moduli}]\label{t-B}
Let \((X,B)\to Y\) be a log Calabi--Yau fibration of dimension \(d\), and let \(S\subset B\) be a horizontal component with coefficient \(\mu_S B=1\). Then
\[
\sdeg(S/Y)
\]
is bounded in terms of \(d\).
\end{theorem}

\begin{example}
Let
\[
X=\Pp^1\times \Pp^1 \to Y=\Pp^1
\]
be projection to the first factor, and let \(B=S\), where \(S\) is a general curve of bidegree \((2,n)\). Then the restriction
\[
f|_S\colon S\to Y
\]
has degree \(2\). Therefore
\[
\sdeg(S/Y)=2.
\]
\end{example}

In relative dimension \(2\), the statement becomes especially concrete.

\begin{fact}
If \((X,B)\to Y\) is a log Calabi--Yau fibration with \(\dim X=2\), and \(S\subset B\) has coefficient \(1\), then
\[
\sdeg(S/Y)\le 2.
\]
\end{fact}

\section{A general boundedness conjecture}

A more general conjecture predicts boundedness under weaker assumptions on the coefficient of \(S\). The following version is slightly more general than that stated in \cite{B-moduli} when $k$ is not algebraically closed.

\begin{conjecture}[Birkar]
Let \((X,B)\to Y\) be a log Calabi--Yau fibration of dimension \(d\) over a field \(k\) of characteristic \(0\), and let \(S\subset B\) be a horizontal component with
\[
\mu_S B\ge t>0.
\]
Then
\[
\sdeg(S/Y)
\]
is bounded in terms of \(d\) and \(t\).
\end{conjecture}

A particularly interesting special case is when \(Y=\Spec k\). Then
\[
\sdeg(S/Y)=\dim_k H^0(S,\OO_S),
\]
so the conjecture predicts a bound on the arithmetic complexity of \(S\). In this form the problem naturally lies at the intersection of algebraic and arithmetic geometry.

\begin{remark}
Of course one can ask whether the conjecture also holds when $S$ is not horizontal but $X\to Y$ is of Fano type. Additionally one can ask if it holds when $k$ has positive characteristic. A mixed characteristic version can also be formulated.
\end{remark}

\begin{remark}
The conjecture applies more generally to closed subsets \(S\subset X\) horizontal over $Y$ where the singularities of the pair are sufficiently bad, i.e. centre of prime divisors with log discrepancy $\le 1-t$ (for example horizontal components of $S$ with $\mu_SB\ge t$ satisfy this condition). This becomes especially interesting when those bad singularities are centered at closed points, since residue field extensions then enter the picture.
\end{remark}

\section{Recent progress}

The following result has recently been obtained.

\begin{theorem}[Birkar--Qu, \cite{BQ-sdeg}]\label{t-BQ}
The conjecture above holds when the base field \(k\) is algebraically closed. More strongly,
\[
\sdeg(S^\nu/Y)
\]
is bounded, where \(S^\nu\) denotes the normalization of \(S\).
\end{theorem}

\begin{remark}
If \(S\) is not horizontal, then \(\sdeg(S/T)\), where \(T\) is the image of \(S\), is not bounded in general without additional assumptions \cite{BQ-sdeg-non-Fano}.
\end{remark}

One surprising consequence is that if \((X,B)\to Y\) and \(S\) satisfy the assumptions of the theorem, and \(F\) is a general fiber of \(X\to Y\), then the restriction \(S|_F\) has only boundedly many irreducible components.

The proofs of both Theorems \ref{t-B} and \ref{t-BQ} rely on sophisticated techniques and results in birational geometry. 

\section{Sketch of the proof strategy}

We now give a brief outline of the proof of Theorem \ref{t-BQ} omitting a lot of details. This is just to get some ideas of how the proof works. For full details, see \cite{BQ-sdeg}. 

\subsection{Horizontal and vertical cases}

Let \(H(d,t)\) denote the horizontal boundedness statement in relative dimension \(d\), and let \(V(d,t)\) denote a corresponding statement for a vertical component \(S\) of \(B\) mapping onto a divisor in \(Y\), in a Fano-type setting.
The idea is to prove both statements in an inductive manner.

\subsection{From \(H(1,t)\) to \(V(1,t)\)}

It is easy to verify $H(1,t)$ as in this case the general fibres of $X\to Y$ are just $\mathbb{P}^1$. One then proves
\[
H(1,t)\Rightarrow V(1,t).
\]
This reduction uses the minimal model program \cite{BCHM} and the theory of complements \cite{B-compl} in subtle ways (somewhat similar to the next step), which allows one to replace the original boundary $B$ by a more controlled one without losing the essential birational information.

\subsection{From \(V(1,t)\) to \(V(d,t)\)}

The next step is to deduce
\[
V(1,t)\Rightarrow V(d,t).
\]
By taking a suitable hyperplane section, one reduces to the case
\[
\dim Y=2,\qquad \dim T=1,
\]
where \(T\) is the image of \(S\) in \(Y\).

Using the minimal model program and complement theory, one can modify the setting and arrange a boundary \(C\) such that:
\begin{itemize}
    \item \((X,C)\) is lc, and \((X,0)\) is \(\epsilon\)-lc over $Y\setminus T$ for some fixed \(\epsilon>0\),
    \item \(K_X+C\sim_\Q 0/Y\),
    \item the coefficients of \(C\) belong to a finite set,
    \item \(S\) supports the whole fiber over \(T\), and 
    \item \(tS\le C\).
\end{itemize}

After further MMP arguments, one reduces to the case when there is a Fano fibration
\[
X\to Z/Y.
\]
If \(\dim(Z/Y)\ge 1\), then one applies the canonical bundle formula together with generalised pairs \cite{BZh} and induction on dimension. 

Thus the key remaining case is \(\dim(Z/Y)=0\), so that \(Z\to Y\) is an isomorphism, perhaps after shrinking $Y$ near the generic point of $T$.

The general fibers of \(X\to Y\) are then \(\epsilon\)-lc Fano varieties, and hence bounded by the BAB theorem \cite{B-BAB}. Using logarithmic geometry \cite{BQ-degen, Q}, one can further reduce to a diagram
\[
\begin{array}{ccc}
X' & \longrightarrow & X'' \\
\downarrow & & \downarrow \\
Y & = & Y
\end{array}
\qquad\text{with a birational map }X\dashrightarrow X'',
\]
such that:
\begin{itemize}
    \item \(X''\to Y\) is a relatively bounded family,
    \item \((X',D')\to Y\) is toroidal, and relatively bounded over an open set, and
    \item $X'\to X''$ is birational.
\end{itemize}

If \(P\) denotes the center of \(S\) on $X'$, then \(P\) becomes an lc center of \((X',D')\). It is therefore enough to bound the number of irreducible components of general fibres of $P\to T$, 
and this can be done using the toroidal structure together with further numerical properties and adjunction.

\subsection{From \(V(d,t)\) to \(H(d,t)\)}

The final step is to prove
\[
V(d,t)\Rightarrow H(d,t).
\]
Again, after MMP and complement theory, one may modify the setting and assume there is a boundary \(C\) such that:
\begin{itemize}
    \item \((X,C)\) is lc and \((X,0)\) is \(\epsilon\)-lc for some fixed $\epsilon>0$,
    \item \(K_X+C\sim_\Q 0/Y\),
    \item \(tS\le C\), and 
    \item there is a Fano fibration \(X\to Z\to Y\).
\end{itemize}

If \(\dim(Z/Y)\ge 1\), then induction applies:
\begin{itemize}
    \item if \(S\) is horizontal over \(Z\), one applies induction to \(X\to Z\),
    \item if \(S\) is vertical over \(Z\), one applies induction both to \(X\to Z\) and to \(Z\to Y\).
\end{itemize}

If instead \(\dim(Z/Y)=0\), then we can assume \(Z\to Y\) is an isomorphism. In this case the general fiber of \(X\to Y\) is an \(\epsilon\)-lc Fano variety, so boundedness follows from BAB. This boundedness is then used to conclude boundedness of
\[
\sdeg(S^\nu/Y).
\]

\begin{remark}
The idea of using toroidal structures has also been applied to establish main results of \cite{B-sing-Fano-fib} and \cite{BQ-degen}.
\end{remark}

{Yau Mathematical Sciences Center, JingZhai, Tsinghua University, Hai Dian District, Beijing, China 100084.\\}
{birkar@tsinghua.edu.cn}

\end{document}